\documentclass[12pt]{article}
\usepackage{amsmath}
\usepackage{amsthm}
\usepackage{amssymb}
\usepackage[latin2]{inputenc}
\usepackage{graphics}

\newtheorem{thm}{\bf Theorem}[section]

\theoremstyle{remark}
\newtheorem{rem}{\bf Remark}[section]

\newtheorem{example}{\bf Example}[section]
\usepackage{epsfig}
\usepackage{color}

\def\a{{\alpha}}

\begin{document}
\title{\Huge Cross--diffusion Modeling in Macroeconomics}
\bigskip
\bigskip
\bigskip
\author{
by\\László Balázsi\footnote{E-mail: Balazsi\_Laszlo@student.ceu.hu} and Krisztina Kiss\footnote{E-mail: kk@math.bme.hu
}\\
Institute of Mathematics, Budapest University of Technology and
Economics,\\H-1111 Budapest, Hungary\\ }
\date{}
\maketitle
\newpage
\begin{abstract}
This paper deals with the stability properties of a closed market, where capital and labour force are acting like a predator--prey system in population--dynamics. The spatial movement of the capital and labour force are taken into account by cross--diffusion effect. First, we are showing two possible ways for modeling  this system in only one country's market (applying a simple functional response and a Holling-type ratio--dependent response as well), examining the conditions of their stability properties. We extend the ratio--dependent model into two countries common market where two kind of cross--diffusion effects are present, and find those additional conditions, whose are necessary for the stability of the global common market besides the stability of each countries local markets. Our four--dimensional model highlights that a hectic movement of the capital toward labour force can cause a Turing instability.

\end{abstract}
\pagestyle{myheadings} \markboth{\centerline {\scriptsize L. Balázsi, K. Kiss
}}{\centerline {\scriptsize  Cross--diffusion modeling in macro--economical model}}


\newpage
\section{Model set up and preliminaries}
\bigskip
\medskip

It is known by Goodwin \cite{Goodvin1} that some economic situation can be described by a Lotka--Volterra predator-prey system. Capital plays the role of prey species while labour force is the predator species. Instead of a simple Lotka--Volterra system it is often used a logistic growth rate of capital \cite{Farkas Kotsis}.
The inhomogeneous distribution of capital and labour force has an important role as it was shown in \cite{Fourier modszer}, \cite{Shaban kozg}. 
The model (1.3) of \cite{Fourier modszer} was modified by cross--diffusion and this effect was modeled by a PDE. In \cite{Shaban kozg} the same model was studies in a patchy environment, namely, the diffusion was modeled by an ODE. It was shown in both cases that cross--diffusion usually causes Turing instability. In our paper we state a different model in which both kind of cross--diffusion effect are present. Basically we are interested in such cases when the diffusion cannot cause instability in our higher dimensional model.

In papers \cite{Fourier modszer}, \cite{Shaban kozg} the kinetic system was modeled by a predator--prey situation as well. We revise the kinetic system according to the following situation.
Let us consider a country, where a firm takes place. It has incentives to expand using its capital stock to invest. For these new firms labour force is required. The total labour force, i.e. the number of those employed and those unemployed, is denoted by $v(t)$ at time $t$. The company supplies free jobs, which number is denoted by $u(t)$ at time $t$ and it is proportional to the invested capital, thus, both labour force and free jobs are attracts each other for sustenance.  We assume that the quantity of free jobs would grow according to the logistic law if there was no labour force available. The labour force decreases the per capita growth rate of free jobs proportional to the quantity of labour force. The factor of proportionality is given by $m>0$. This is the rate by which the labour force is filling in the free jobs. This leads to the following differential equation:
\begin{equation}\label{h:01}
\dot{u}(t)=u(t) r \left( 1-\frac{u(t)}{K}\right)-m u(t) v(t)\\
\end{equation}
where $r$ is the natural per capita growth rate of free jobs and $K>0$ is its theoretical eventual maximum.  Equation (\ref{h:01}) is exactly the same as (1.1) in \cite{Fourier modszer}.
We assume that the per capita growth rate of labour force is increasing proportionally with respect to the free jobs and we consider the case when the factor of proportionality is equal to $m$. Moreover we suppose that in the absence of free jobs, the labour without revenue is unable to maintain itself, leading to its declension. The death rate of labour force is denoted by $d>0$.
This leads to the differential equation:
\begin{equation}\label{h:02}
\dot{v}(t)=m u(t) v(t)-d
 v(t).
\end{equation}
This equation is different from the equation (1.2) given in \cite{Fourier modszer}.
In population dynamics the function $m v(t)$ is called functional response. It is often fails to be linear. Later we also show such model in which the functional response is more complicated. But in this simple case the following system of differential equations characterizes the situation:
\begin{equation}\label{h:0}\left.\begin{array}{lll}
\dot{u}(t)&=&u(t) r \left( 1-\frac{u(t)}{K}\right)-m u(t) v(t)\\
\dot{v}(t)&=&m u(t) v(t)-d v(t)
\end{array}\right\}
\end{equation}
Comparing model (\ref{h:0}) to the one given in \cite{Fourier modszer}, we can see that both model is a predator--prey situation where capital is represented by prey species while labour force is by predator species. The only difference between them that in the model given in \cite{Fourier modszer}, both prey and predator can survive without the other, while in (\ref{h:0})  the predator dies out in absence of prey. Moreover, model (\ref{h:0}) is the same as model (2.5) without delay given in \cite{Farkas Kotsis}.
As in \cite{Fourier modszer} we modify model (\ref{h:0}) by taking into account that the economy of the country is not concentrated in a point but is disturbed in a spatial bounded domain.  It is assumed that capital and labour force is moving around freely in the country and that the capital investment and the conditions of living are the same everywhere. Besides, we assume that the economy is closed, meaning that there is no in-- and outflow of capital and labour force at the boundary.
In order to be able to perform explicit calculations we consider the spatial domain, abstractly, as one dimensional and identify it with the interval $[0,L]\subset\mathbb{R}$, $L>0$. From now on, $u(t,x)$ and $v(t,x)$ represent the densities of free jobs and labour forces, respectively, at time $t\geq 0$ at the point $x\in [0,L]$ of the domain.
We assume that the movement of capital and labour force follow Fick's law, i.e. both substances flow away from places where their densities are high, towards places where their densities are low. Furthermore, both substances are influenced by each other: larger quantity of free jobs attracts labour force rather then smaller, and new firms are more likely to be founded in places, where the density of labour force is higher. The velocities and the direction of flowing can be represented by diffusion coefficients and their signs, contained by the diffusion matrix:
\begin{equation}\label{kor:1}
\mathbf{D}=\left(
\begin{array}{ccc}
d_{11} & -d_{12}\\
-d_{21} & d_{22}
\end{array}\right),
\end{equation}
where $d_{11}$ and $d_{22}$ denote the self--diffusion of capital and labour force respectively, $d_{12}$ the velocity of capital--flowing towards labour force, and $d_{21}$ shows the largeness of the migration of labour force towards free jobs and $d_{ij}>0$, $i,j=1,2$, denote the self--diffusion effect (e.g. \cite{Keresztdiff}). In order to obtain a well posed problem, necessary to have $\det\mathbf{D}>0$. (The same diffusion effect was considered in \cite{Fourier modszer}.)
Now we can set up our first reaction--diffusion model:
\begin{equation}\label{h:1}\left.\begin{array}{lll}
u'_t(t,x)&=&u(t,x) r \left( 1-\frac{u(t,x)}{K}\right)-m u(t,x) v(t,x)+d_{11} u''_{xx}(t,x)-d_{12} v''_{xx}(t,x)\\
v'_t(t,x)&=&m u(t,x) v(t,x)-d v(t,x)-d_{21} u''_{xx}(t,x)+d_{22} v''_{xx}(t,x)
\end{array}\right\},\end{equation}
\begin{equation}\label{h:1a}
u^{'}_x(t,0)=u^{'}_x(t,L)=0=v^{'}_x(t,0)=v^{'}_x(t,L),
\end{equation}
where boundary condition (\ref{h:1a}) presents the closeness of the market.
The initial conditions are given by
\begin{equation}\label{h:1b}
u(0,x)=U(x), \;\; v(0,x)=V(x), \;\; x\in[0,L], \;\;U,V\in \mathbb{C}^2([0,L], \mathbb{R}).
\end{equation}
It will be  shown that at certain values of $d_{ij}$, $i,j=1,2$ a Turing bifurcation may occur as it often happens in similar situation.

This paper is organized as follows. In section \ref{simple} we are giving an overview of the stability behaviour of model (\ref{h:1}). It will be shown that this model has some disadvantageous properties. In section \ref{ratio}  model (\ref{h:1}) will be modified by a so-called ratio--dependent functional response and the effect of the changing of the diffusion coefficients will be studied to the stability behaviour. In section \ref{patch} we are extending the basic model to a patchy environment illustrating the outer movement of capital and labour force between two countries without and with inner migration. We state and study a four--dimensional PDE model which is capable to describe the inner and outer movements of the market at the same time.
The model (\ref{h:1rd}) in section \ref{ratio} is a well known model and model (\ref{p:2}) in section \ref{patch} is also known. We repeat here partly the known stability investigation of them and we emphasize some stability conditions by other aspects in order to understand better our new model (\ref{pi:1hat}) in section \ref{patch}.
 In section \ref{sum}
we are summarizing the main conclusions of the study.

\section{Stability of the simple market}\label{simple}

The kinetic system of (\ref{h:1}) has three equilibria, namely:
$$E_0=(0,0); \;\;\; E_1=(K,0); \;\;\; \bar{E}=(\bar{u},\bar{v}),$$
where
\begin{equation}\label{eq}
(\bar{u},\bar{v})=(\frac{d}{m},\frac{-d r+K m r}{K m^{2}})
\end{equation}
and $\bar{v}>0$ if $K>\frac{d}{m}$.
Each of them can be examined, however we have special interest in $\bar{E}$ (given by (\ref{eq})), because this is the only one which is not placed on the boundary of the plane. The equilibrium of the kinetic system
is also a constant solution of (\ref{h:1}). The interaction matrix of the kinetic system linearized at $\bar{E}$ is:
\begin{equation*}
\mathbf{J_{\bar{E}}}=\mathbf{A}=\left(
\begin{array}{cc}
\frac{-d r}{Km} & -d \\
r \left(1-\frac{d}{Km}\right) & 0
\end{array}\right).
\end{equation*}
It is easy to see, that this is in fact stable, if
\begin{equation}\label{h:2}K>\frac{d}{m}.\end{equation}
Thus, the following theorem holds:
\begin{thm}\label{T2.1}
Suppose, that condition (\ref{h:2}) is satisfied; then the equilibrium $\bar{E}=(\bar{u},\bar{v})$ of the kinetic system of (\ref{h:1}) is asymptotically stable.
\end{thm}

Theorem \ref{T2.1} shows that system (\ref{h:1}) is a too simple model because if a positive equilibrium exists then it is necessarily stable. (We note here that if condition (\ref{h:2}) holds then matrix $\mathbf{J_{\bar{E}}}$ is signstable e.g. \cite{vdDriessche}.)

Now we are linearizing the reaction--diffusion system (\ref{h:1}) around the constant solution $(\bar{u},\bar{v})$, using the notation:
\begin{equation*}p(t,x)=u(t,x)-\bar{u};\qquad q(t,x)=v(t,x)-\bar{v}.\end{equation*}
We get
\begin{equation}\label{h:3}\left.\begin{array}{lll}
p'_t(t,x)&=&\left(r-\frac{2r\bar{u}}{K}-m\bar{v}\right) p(t,x)-(m\bar{u}) q(t,x)+d_{11} p''_{xx}(t,x)-d_{12} q''_{xx}(t,x)\\
q'_t(t,x)&=&(m\bar{v}) p(t,x)+(m\bar{u}-d) q(t,x)-d_{21} p''_{xx}(t,x)+d_{22} q''_{xx}(t,x)
\end{array}\right\}.\end{equation}
The boundary conditions are obtained:
\begin{equation}\label{h:3a}
p^{'}_x(t,0)=p^{'}_x(t,L)=0=q^{'}_x(t,0)=q^{'}_x(t,L).
\end{equation}
Solving (\ref{h:3}) by Fourier--method, it is easy to see (e. g. \cite{FM BIOL}, \cite{Fourier modszer}) that the stability behaviour depends on the stability of the following matrices:
\begin{equation*}
\mathbf{A}-\lambda_k\mathbf{D}=\left[\begin{array}{cc}\frac{-dr}{Km}-\lambda_k d_{11} & -d+\lambda_k d_{12} \\
r \left(1-\frac{d}{Km}\right)+\lambda_k d_{21} & -\lambda_k d_{22}\end{array}\right], \;\; \lambda_k=(\frac{k\pi}{L})^2, \;\; k=0,1,2,\dots.
\end{equation*}
Matrix $\mathbf{A}-\lambda_k\mathbf{D}$ is stable if the following conditions are satisfied:
\begin{equation*}
Tr(\mathbf{A}-\mathbf{D}\lambda_k)=Tr(\mathbf{A})-\lambda_k Tr(\mathbf{D})=-\underbrace{\frac{dr}{Km}}_{>0}-\underbrace{\lambda_k (d_{11}+d_{22})}_{>0}<0
\end{equation*}
\begin{equation*}
\det\left(\mathbf{A}-\mathbf{D}\lambda_k\right)=\underbrace{\det\mathbf{A}}_{>0}+\lambda^2_k \underbrace{\det\mathbf{D}}_{>0}-\lambda_k(-\frac{dr}{Km} d_{22}-d d_{21}+r\left(1-\frac{d}{Km}\right) d_{12})>0.
\end{equation*}
Since $\det\mathbf{A}>0$ because of (\ref{h:2}) and we assumed initially, that $\det\mathbf{D}>0$, the sufficient condition to get a positive determinant for all $\lambda_k$ is
\begin{equation*}
-\frac{dr}{Km} d_{22}-d d_{21}+r\left(1-\frac{d}{Km}\right) d_{12}<0.
\end{equation*}
(This is the case when $\det\left(\mathbf{A}-\mathbf{D}\lambda_k\right)$ is a stable polynomial.)
Expressing $d_{12}$ from the above condition, we obtain:
\begin{equation}\label{h:5}
d_{12}<\frac{\frac{dr}{Km} d_{22}+d d_{21}}{r\left(1-\frac{d}{Km}\right)}.
\end{equation}
Thus, we have arrived at:
\begin{thm}
 If conditions (\ref{h:2}) and (\ref{h:5}) hold then the constant solution $(0,0)$ of the linear problem (\ref{h:3})--(\ref{h:3a}) is asymptotically stable.
\end{thm}
\begin{rem}
According to \cite{Casten Holland} (or to \cite{smoller}), the asymptotic stability of the constant solution $(\bar{u},\bar{v})$ of system (\ref{h:1})--(\ref{h:1a}) is implied.
\end{rem}
Increasing the diffusion coefficient $d_{12}$, exceeding the critical value:
\begin{equation*}
d_{12_{crit}}=\frac{\frac{dr}{Km} d_{22}+d d_{21}}{r\left(1-\frac{d}{Km}\right)},
\end{equation*}
where condition (\ref{h:5}) does not hold any more, equilibrium $(\bar{u},\bar{v})$ of system (\ref{h:1})--(\ref{h:1a}) may undergo a Turing--bifurcation \cite{turinst} and became unstable. More precisely, if $d_{12}>d_{12_{crit}}$ then there exists $\lambda_k$ for which  $\det\left(\mathbf{A}-\mathbf{D}\lambda_k\right)<0$. 
 Examining the sign of the diffusion coefficients in condition (\ref{h:5}), it is clear that higher values of $d_{21}$ and $d_{22}$ are supporting the stability of the system, while $d_{12}$ has a destabilizing effect as it increases. Hence, our calculations lead us to the following conclusion: rapid movement of labour force towards free jobs and fast moving of labour force to areas where its densities are low from places where its densities are high are advantageous for the stability of the market, while the hectic movement of capital stock (free jobs) towards places where the labour force density is higher acts negatively, even has the ability to make the market unstable. Interesting fact that diffusion coefficient $d_{11}$ does not have an influence in stability in this case. Namely, the movement of free jobs from the places where its density is high towards the places where its density is low has no importance. This latter is also a disadvantageous characteristic of this model. Thus, we turn to study our model using another functional response.

\section{Model with Holling--type ratio--dependent response}\label{ratio}

In population dynamics a Holling-type ratio--dependent functional response is used in that cases when the competition for food is very sharp, see for example in \cite{AkcakayaArditiGinzburg95}, \cite{ArditiBerryman91}, \cite{AriditiGrinzburg89}, \cite{CosnerDeAngelisAultOlson99}, \cite{Shaban ratio}, \cite{KK-KS}, \cite{KK-Lizana}. It is reasonable to apply this kind of functional response when the labour force need to search for free jobs. According to this, the following model is relevant when the competition for free jobs is intensive:
\begin{equation}\label{h:1rd}\left.\begin{array}{lll}
u'_t(t,x)&=&u(t,x) r \left( 1-\frac{u(t,x)}{K}\right)-m \frac{u(t,x)v(t,x)}{av(t,x)+u(t,x)}+d_{11} u''_{xx}(t,x)-d_{12} v''_{xx}(t,x)\\
v'_t(t,x)&=&m \frac{u(t,x)v(t,x)}{av(t,x)+u(t,x)}-d v(t,x)-d_{21} u''_{xx}(t,x)+d_{22} v''_{xx}(t,x)
\end{array}\right\},\end{equation}
with the boundary conditions (\ref{h:1a}) and the initial conditions (\ref{h:1b}).
Parameter $a>0$ is called half--saturation constant (e. g. \cite{KK-KS}) and $m>0$ is the maximum growth rate of labour force. Of course the model is relevant if the maximal growth rate of the labour force is larger than its death rate, namely
\begin{equation}\label{h:3rd}
m-d>0.
\end{equation}
\\
The equilibria of the kinetic system of (\ref{h:1rd}) are:
$$E_1=(K,0); \;\;\; \bar{E_r}=(\bar{u},\bar{v}),$$
where
\begin{equation}\label{eqrd}
(\bar{u},\bar{v})=(\frac{K(d-m+ar)}{ar},\frac{K (m-d) (d-m+ar)}{a^2dr})
\end{equation}
\begin{rem}
Equilibrium $E_0=(0,0)$ is missing, because certain denominators would be equal to zero. This problem could be solved by continuous extension of the model (e. g. \cite{KK-KS}), however we are interested in $\bar{E_r}$.
\end{rem}
Necessary to have the equilibrium $\bar{E_r}$ positive, which requires the following condition:
\begin{equation}\label{h:2rd}
r>\frac{m - d}{a}.
\end{equation}
This means that the natural per capita growth rate of free jobs (capital) has to be higher than the ratio of the maximal growth rate and half saturation constant of the labour force. This is natural because high growth rate of labour force (relatively to the half saturation constant) would result a too fast falling in free jobs. \\
Following the method given in section \ref{simple},
linearizing (\ref{h:1rd}) in $\bar{E_r}$ (given by (\ref{eqrd})) we obtain the interaction matrix in the following form:
\begin{equation*}
\mathbf{J_{\bar{E}_r}}=\mathbf{A_r}=\left(
\begin{array}{cc}
\frac{m^2-d^2}{ma}-r & -\frac{d^2}{m}\\
\frac{(d-m)^2}{am} & -\frac{d(m-d)}{m}
\end{array}\right).
\end{equation*}
This is stable if $\det\mathbf{J_{\bar{E}_r}}>0$ and $Tr(\mathbf{J_{\bar{E}_r}})<0$, namely
\begin{equation*}
\det\mathbf{J_{\bar{E}_r}}=\frac{d (m-d) (d-m+ar)}{am}>0,
\end{equation*}
which holds if and only if conditions (\ref{h:2rd}) and (\ref{h:3rd}) hold.
\begin{equation*}\begin{split}
Tr(\mathbf{J_{\bar{E}_r}})=\frac{(a-1) d^2-adm+m (m-ar)}{am}\underbrace{<}_{by\;(\ref{h:2rd})}\frac{(a-1) d^2-adm+md}{am}=\\=\frac{(a-1) (d^2-dm)}{am}=\frac{(a-1) d \overbrace{(d-m)}^{<0\;by \;(\ref{h:3rd})}}{am}<0
\end{split}\end{equation*}
holds if
\begin{equation}\label{h:4rd}
a>1.
\end{equation}
In population dynamics the low value of parameter $a$ refers to a predator species which tries to ensure its survival by having a relatively low growth rate and consuming less. It is a so--called K--strategist (e.g.\cite{FM BIOL}). In our case the meaning of (\ref{h:4rd}) will be that we apart from low willingness to work, low desire to work, low interest to work or discharged worker.
Summing up the conditions we have obtained:
\begin{thm}\label{theoratiokin}
Suppose that conditions (\ref{h:3rd}), (\ref{h:2rd}) and (\ref{h:4rd}) hold; thus equilibrium $\bar{E_r}$ (given by (\ref{eqrd})) of the kinetic system of (\ref{h:1rd}) is asymptotically stable.
\end{thm}


 Now, just as in section \ref{simple}, we solve the linearized model of (\ref{h:1rd})--(\ref{h:1a}) by Fourier--method. We get that the stability behaviour of the spatially constant solution depends on the stability of the following matrices:
\begin{equation*}
\mathbf{A_r}-\mathbf{D}\lambda_k=\left[
\begin{array}{cc}
\frac{m^2-d^2}{ma}-r-\lambda_kd_{11} & -\frac{d^2}{m}+\lambda_kd_{12}\\
\frac{(d-m)^2}{am}+\lambda_kd_{21} & -\frac{d(m-d)}{m}-\lambda_kd_{22}
\end{array}\right], \;\; \lambda_k=(\frac{k\pi}{L})^2, \;\; k=0,1,2,\dots.
\end{equation*}
In order to be this matrix stable the following conditions need to be satisfied for all $\lambda_k$:
\begin{equation*}
Tr(\mathbf{A_r}-\mathbf{D}\lambda_k)<0; \qquad \det\left(\mathbf{A_r}-\mathbf{D}\lambda_k\right)>0.
\end{equation*}
From condition (\ref{h:4rd}) we know that $Tr(\mathbf{A_r})<0$, using this:
\begin{equation*}
Tr(\mathbf{A_r}-\mathbf{D}\lambda_k)=Tr(\mathbf{A_r})-\lambda_k \underbrace{(d_{11}+d_{22})}_{>0}<0.
\end{equation*}
Let us consider
\begin{equation*}
\det\left(\mathbf{A_r}-\mathbf{D}\lambda_k\right)=\det\mathbf{A_r}+\lambda^2_k \det\mathbf{D}-\lambda_k (a_{r_{11}} d_{22}+a_{r_{22}} d_{11}+a_{r_{12}} d_{21}+a_{r_{21}} d_{12})>0,
\end{equation*}
where $A_r=[a_{r_{ij}}]$,
$\det\mathbf{D}>0$, and if conditions (\ref{h:3rd}) and (\ref{h:2rd}) hold, we have $\det\mathbf{A_r}>0$ as well, thus, a sufficient condition to have a positive determinant is
\begin{equation}\label{ar-es}
a_{r_{11}} d_{22}+a_{r_{22}} d_{11}+a_{r_{12}} d_{21}+a_{r_{21}} d_{12}<0.
\end{equation}
In this case $\det\left(\mathbf{A_r}-\mathbf{D}\lambda_k\right)$ is a stable polynomial.
Substituting the given entries of $A_r$ and expressing $d_{12}$:
\begin{equation}\label{h:7rd}
d_{12}<\frac{ad}{m-d} d_{11}+\frac{ad^2}{(m-d)^2} d_{21}+\left(-\frac{m+d}{m-d}+\frac{amr}{(m-d)^2}\right) d_{22}.
\end{equation}
Applying again the theorem of Casten--Holland  \cite{Casten Holland}  we get:
\begin{thm}
Suppose that conditions (\ref{h:3rd}), (\ref{h:2rd}), (\ref{h:4rd}), (\ref{h:7rd}) hold then the constant solution $E_r$ of system (\ref{h:1rd})--(\ref{h:1a}) is asymptotically stable.
\end{thm}
\begin{rem}
The coefficient of $d_{22}$ takes positive values in (\ref{h:7rd}) if
\begin{equation}\label{plus}
\frac{amr}{(m-d)}>m+d,
\end{equation}
or
\begin{equation}\label{plusmas}
r>\frac{m-d}{a}(1+\frac{d}{m})
\end{equation}
holds.
From this the positivity of $E_r$ follows. In this case $a_{r_{11}}<0$. This means that we are outside of the so--called Allee--effect zone of free jobs. (See the meaning of this zone in \cite{FM BIOL}, \cite{KK-KS} in connection with a prey species.) Outside of the Allee--effect zone the number of free jobs are not too low.
\end{rem}
Increasing diffusion coefficient $d_{12}$, exceeding a critical value, where condition (\ref{h:7rd}) does not hold any more:
\begin{equation*}
d_{12_{crit}}=\frac{ad}{m-d} d_{11}+\frac{ad^2}{(m-d)^2} d_{21}+\left(-\frac{d+m}{m-d}+\frac{amr}{(m-d)^2}\right) d_{22},
\end{equation*}
system (\ref{h:1rd})--(\ref{h:1a}) may undergo a Turing--bifurcation and further increase of $d_{12}$ causes instability. (Similarly as in section \ref{simple}.)
 Examining the diffusion coefficients in condition (\ref{h:7rd}), it is obvious that higher values of $d_{11}$, $d_{22}$ and  $d_{21}$ (if (\ref{plus}) holds), are supporting the stability of the market, while an increase in coefficient $d_{12}$ affects negatively on it. Namely,  fast movement of firms towards rarely industried areas supports the stability of the market,
 so as rapid migration of labour force towards places where its densities are low in case of greater unemployment and where free jobs are abundant. Nevertheless, hectic flow of the capital towards places where labour force is dense acts negatively and also able to violate the stability of the system.

\section{The patchy model}\label{patch}

 In this section we study the cross--diffusion effect between two countries.  They will be modeled by patches e.g. \cite{Shaban Farkas PP, Shaban Farkas Comp, Shaban kozg}. We apply these methods for our cases as follows: Let us consider two countries as two patches, providing the following: the inner migration remains still inside each patch, but a new movement of capital and labour force has taken into account between the two countries. We are curious about the conditions of the stability of the common market, providing that each market is also stable, and about how the diffusion affects on the obtained stability. The Holling--type ratio--dependent model introduced in section \ref{ratio} will be used.

\subsection{Without inner migration}\label{withoutim}

First we consider the patchy model as two points, thus we neglect the inner migration of the two substances. The model can be written as follows:
\begin{equation}\label{p:2}\left.\begin{array}{lll}
\dot{u}(1,t)&=&f(u(1,t),v(1,t))+\delta_1 (\rho_1(v(2,t))u(2,t)-\rho_1(v(1,t))u(1,t))\\
\dot{v}(1,t)&=&g(u(1,t),v(1,t))+\delta_2 (\rho_2(u(2,t))v(2,t)-\rho_2(u(1,t))v(1,t))\\
\dot{u}(2,t)&=&f(u(2,t),v(2,t))+\delta_1 (\rho_1(v(1,t))u(1,t)-\rho_1(v(2,t))u(2,t))\\
\dot{v}(2,t)&=&g(u(2,t),v(2,t))+\delta_2 (\rho_2(u(1,t))v(1,t)-\rho_2(u(2,t))v(2,t))
\end{array}\right\},\end{equation}
where numbers $1,2$ in the arguments refer to the number of patch, thus, those mean the number of the country, while $f$ and $g$ are the same functions as in the kinetic system of the ratio--dependent model (\ref{h:1rd}), namely:
\begin{equation}\label{p:3}\left.\begin{array}{lll}
f(u,v)&=&ru\left(1-\frac{u}{K}\right)-m \frac{uv}{av+u}\\
g(u,v)&=&m \frac{uv}{av+u}-dv
\end{array}\right\},\end{equation}
$\rho_1, \rho_2 \in \mathbb{C}^1$ are the migration functions. Similar model was studied in \cite{Shaban kozg} but now the reaction term (\ref{p:3}) is different. These latter functions describe the diffusion of capital and labour respectively between the two countries (patches), $\delta_i$ is the diffusion coefficient, denotes the velocity of this movement. Functions $\rho_1, \rho_2>0$, and  $\rho'_1, \rho'_2<0$. (Here $\rho'_1(v(i,t))=\frac{d \rho_1(v(i,t))}{d v(i,t)}, \;\;\rho'_2(u(i,t))=\frac{d \rho_2(u(i,t))}{d u(i,t)}, \;\; i=1,2$.) These positive  decreasing functions expressing the facts that the migration of the labour is (free jobs are) higher into a country from the other if the number of the free jobs  are (labour is) less in the other country. These functions model the cross--diffusion effect between the countries. We say that the cross--diffusion is strong if $|\rho_i'|$ is large. If $\rho_i\equiv 1$ then we have mere self--diffusion.
Since problem (\ref{p:2}) has the same kinetic system as problem (\ref{h:1rd})--(\ref{h:1a}) if we restrict it to a single country, all the equilibria and stability conditions will remain unchanged. Let us denote the interaction matrix of system (\ref{p:2}) at $\bar{E_p}=(\bar{u},\bar{v},\bar{u},\bar{v})$ by $\mathbf{A_p}$ (where $\bar{u},\bar{v}$ are given by (\ref{eqrd})), we get:
\begin{equation*}
\mathbf{A_p}=\left(\begin{array}{cccc} a_r{_{11}}& a_r{_{12}}&0&0\\ a_r{_{21}}& a_r{_{22}}&0&0\\0&0&a_r{_{11}}& a_r{_{12}}\\0&0&a_{21}& a_r{_{22}}\end{array}\right)=\left(\begin{array}{cccc}\frac{m^2-d^2}{ma}-r& -\frac{d^2}{m}&0&0\\ \frac{(d-m)^2}{am}& -\frac{d(m-d)}{m}&0&0\\0&0&\frac{m^2-d^2}{ma}-r& -\frac{d^2}{m}\\0&0&\frac{(d-m)^2}{am}& -\frac{d(m-d)}{m}\end{array}\right).
\end{equation*}
\begin{thm}
Suppose that conditions (\ref{h:3rd}), (\ref{h:2rd}), (\ref{h:4rd}) hold (for both countries), thus, the equilibrium $\bar{E_r}$ is asymptotically stable in each markets, which implies the asymptotic stability of equilibrium $\bar{E_p}=(\bar{u},\bar{v},\bar{u},\bar{v})$ of the kinetic system of (\ref{p:2}) (where $(\bar{u},\bar{v})$ is given by (\ref{eqrd})).
\end{thm}
\par Now denote the coefficient--matrix of the linearized reaction--diffusion system  (\ref{p:2}) at $\bar{E_p}$  by $(\mathbf{A_p}+\mathbf{\Gamma})$, where
\begin{equation*}\mathbf{\Gamma}=\left(\begin{array}{cccc}
-\delta_1\rho_1(\bar{v})&-\delta_1\rho'_1(\bar{v})\bar{u}&\delta_1\rho_1(\bar{v})&\delta_1\rho'_1(\bar{v})\bar{u}\\
-\delta_2\rho'_2(\bar{u})\bar{v}&-\delta_2\rho_2(\bar{u})&\delta_2\rho'_2(\bar{u})\bar{v}&\delta_2\rho_2(\bar{u})\\
\delta_1\rho_1(\bar{v})&\delta_1\rho'_1(\bar{v})\bar{u}&-\delta_1\rho_1(\bar{v})&-\delta_1\rho'_1(\bar{v})\bar{u}\\
\delta_2\rho'_2(\bar{u})\bar{v}&\delta_2\rho_2(\bar{u})&-\delta_2\rho'_2(\bar{u})\bar{v}&-\delta_2\rho_2(\bar{u}).
\end{array}\right).\end{equation*}
The characteristic polynomial can be determined by row- and column operations:\\
$\det\left(\mathbf{A_p}+\mathbf{\Gamma}-\mu\mathbf{I}\right)$
\begin{equation*}\begin{split}
=
\left|\begin{array}{cccc}
a_{r_{11}}-\delta_1\rho_1(\bar{v})-\mu & a_{r_{12}}-\delta_1\rho'_1(\bar{v})\bar{u} & \delta_1\rho_1(\bar{v}) & \delta_1\rho'_1(\bar{v})\bar{u} \\
a_{r_{21}}-\delta_2\rho'_2(\bar{u})\bar{v} & a_{r_{22}}-\delta_2\rho_2(\bar{u})-\mu & \delta_2\rho'_2(\bar{u})\bar{v} & \delta_2\rho_2(\bar{u}) \\
\delta_1\rho_1(\bar{v}) & \delta_1\rho'_1(\bar{v})\bar{u} & a_{r_{11}}-\delta_1\rho_1(\bar{v})-\mu & a_{r_{12}}-\delta_1\rho'_1(\bar{v})\bar{u} \\
\delta_2\rho'_2(\bar{u})\bar{v} & \delta_2\rho_2(\bar{u}) & a_{r_{21}}-\delta_2\rho'_2(\bar{u})\bar{v} & a_{r_{22}}-\delta_2\rho_2(\bar{u})-\mu \\
\end{array}\right|\\
\underbrace{=}_{c1+c3;\;\;c2+c4}\left|\begin{array}{cccc}
a_{r_{11}}-\mu & a_{r_{12}} & \delta_1\rho_1(\bar{v}) & \delta_1\rho'_1(\bar{v})\bar{u} \\
a_{r_{21}} & a_{r_{22}}-\mu & \delta_2\rho'_2(\bar{u})\bar{v} & \delta_2\rho_2(\bar{u}) \\
a_{r_{11}}-\mu & a_{r_{12}} & a_{r_{11}}-\delta_1\rho_1(\bar{v})-\mu & a_{r_{12}}-\delta_1\rho'_1(\bar{v})\bar{u} \\
a_{r_{21}} & a_{r_{22}}-\mu & a_{r_{21}}-\delta_2\rho'_2(\bar{u})\bar{v} & a_{r_{22}}-\delta_2\rho_2(\bar{u})-\mu \\
\end{array}\right|\\\underbrace{=}_{r3-r1;\;\;r4-r2}\left|\begin{array}{cccc}
a_{r_{11}}-\mu & a_{r_{12}} & \delta_1\rho_1(\bar{v}) & \delta_1\rho'_1(\bar{v})\bar{u} \\
a_{r_{21}} & a_{r_{22}}-\mu & \delta_2\rho'_2(\bar{u})\bar{v} & \delta_2\rho_2(\bar{u}) \\
0 & 0 & a_{r_{11}}-2\delta_1\rho_1(\bar{v})-\mu & a_{r_{12}}-2\delta_1\rho'_1(\bar{v})\bar{u} \\
0 & 0 & a_{r_{21}}-2\delta_2\rho'_2(\bar{u})\bar{v} & a_{22}-2\delta_2\rho_2(\bar{u})-\mu \\
\end{array}\right|\\=\det\left[\begin{array}{cc}a_{r_{11}}-\mu & a_{r_{12}}\\a_{r_{21}} & a_{r_{22}}-\mu\end{array}\right]
  \det\left[\begin{array}{cc} a_{r_{11}}-2\delta_1\rho_1(\bar{v})-\mu & a_{r_{12}}-2\delta_1\rho'_1(\bar{v})\bar{u}\\a_{r_{21}}-2\delta_2\rho'_2(\bar{u})\bar{v} & a_{22}-2\delta_2\rho_2(\bar{u})-\mu\end{array}\right].
\end{split}
\end{equation*}
The first determinant in the previous product is the characteristic polynomial of the kinetic system of (\ref{h:1rd}), thus, according to theorem \ref{theoratiokin} it is stable if conditions (\ref{h:3rd}), (\ref{h:2rd}) and (\ref{h:4rd}) are satisfied. Consequently we have to determine the cases where the second determinant, meaning a polynomial, will be stable. Let us denote the matrix by $B=[b_{ij}]$  whose characteristic polynomial is given by
 \begin{equation}\label{p:4}
\det\left[\begin{array}{cc} a_{r_{11}}-2\delta_1\rho_1(\bar{v})-\mu & a_{r_{12}}-2\delta_1\rho'_1(\bar{v})\bar{u}\\
a_{r_{21}}-2\delta_2\rho'_2(\bar{u})\bar{v} & a_{r_{22}}-2\delta_2\rho_2(\bar{u})-\mu
\end{array}\right].\end{equation}
If conditions (\ref{h:3rd}), (\ref{h:2rd}), (\ref{h:4rd}), (\ref{plusmas}), $\rho'_{i}<0, \;\;\rho_{i}>0, \;\; i=1,2$ and
\begin{equation}\label{sign}
a_{r_{12}}-2\delta_1\rho'_1(\bar{v})\bar{u}<0
\end{equation}
are satisfied then $B$ is sign--stable e.g. \cite{vdDriessche} .
Substituting $a_{r_{11}}$ and expressing $\delta_1$ from (\ref{sign}) we obtain:
\begin{equation}\label{p:5}
\delta_1<-\frac{d^2}{2\rho'_1(\bar{v})\bar{u}m}.
\end{equation}
 Parameter $\delta_1$ shows the velocity of the capital (free jobs) movement between the countries. According to (\ref{p:5}), $\delta_1$ has to be  under a limit, namely too fast flow of capital is not allowed.
With this consideration, we have arrived at:
\begin{thm}
Suppose that conditions (\ref{h:3rd}), (\ref{h:2rd}), (\ref{h:4rd}), (\ref{plusmas}), $\rho'_{i}<0, \;\;\rho_{i}>0, \;\; i=1,2$ hold for each patch and (\ref{p:5}) also holds for system (\ref{p:2}) then equilibrium $\bar{E_p}=(\bar{u},\bar{v},\bar{u},\bar{v})$  of system (\ref{p:2}) (given by (\ref{eqrd})) is asymptotically stable.
\end{thm}
This theorem highlights, that for the stability of the common market, relatively slow foreign firm--movement and a significant number of free jobs are required besides the original stability conditions for each market individually. Consequently, the hectic cross--country flow of the capital has to be avoided again.\\
If this particular capital investment is fast enough i.e.:
\begin{equation*}
\delta_1\geq-\frac{d^2}{2\rho'_1(\bar{v})\bar{u}m},
\end{equation*}
then the common market may lose its stability.

The stability of (\ref{p:4}) can be ensured by the following way, as well.
The entries in the main diagonal in $B$ are negative, thus, the determinant of $B$ has to be positive in order to get stability:
$\det B = \det A_r -2(\delta_1\rho_1(\bar{v})a_{r_{22}}+\delta_2\rho_2(\bar{u})a_{r_{11}}-\delta_1\rho'_1(\bar{v})\bar{u}a_{r_{21}}-\delta_2\rho'_2(\bar{u})\bar{v}a_{r_{12}})+
4\delta_1\delta_2(\rho_1(\bar{v})\rho_2(\bar{u})-\rho'_1(\bar{v})\bar{u}\rho'_2(\bar{u})\bar{v})$.
If (\ref{h:3rd}) and (\ref{plusmas}) are satisfied, $A_r$ is signstable, thus, $\det A_r>0$.
If
\begin{equation}\label{1.feltetel}
1>\frac{\rho'_1(\bar{v})\bar{u}}{\rho_1(\bar{v})}\frac{\rho'_2(\bar{u})\bar{v}}{\rho_2(\bar{u})}
\end{equation}
then $(\rho_1(\bar{v})\rho_2(\bar{u})-\rho'_1(\bar{v})\bar{u}\rho'_2(\bar{u})\bar{v})>0$.
If
\begin{equation}\label{2.feltetel}
\frac{a_{r_{22}}}{a_{r_{21}}}<\frac{\rho'_1(\bar{v})\bar{u}}{\rho_1(\bar{v})}
\end{equation}
then $(\delta_1\rho_1(\bar{v})a_{r_{22}}+\delta_2\rho_2(\bar{u})a_{r_{11}}-\delta_1\rho'_1(\bar{v})\bar{u}a_{r_{21}}-\delta_2\rho'_2(\bar{u})\bar{v}a_{r_{12}})<0$.
After short calculations we obtain (\ref{2.feltetel}) in the following form:
\begin{equation}\label{2.feltetelujalak}
-\frac{1}{\bar{v}}<\frac{\rho'_1(\bar{v})}{\rho_1(\bar{v})}.
\end{equation}
Now we can summarize our result in the following theorem:
\begin{thm}\label{TH4.3}
Suppose that conditions (\ref{h:3rd}), (\ref{plusmas}), (\ref{1.feltetel}), (\ref{2.feltetelujalak}), $\rho'_{i}<0, \;\;\rho_{i}>0, \;\; i=1,2$ hold for each patch then equilibrium $\bar{E_p}=(\bar{u},\bar{v},\bar{u},\bar{v})$  of system (\ref{p:2}) (given by (\ref{eqrd})) is asymptotically stable.
\end{thm}
\begin{example}
In the case of $\rho_1({v})=\frac{\alpha+v}{1+v}$, $\rho_2({u})=\frac{\beta+u}{1+u}$ $\alpha, \beta>1$ conditions (\ref{1.feltetel}), (\ref{2.feltetelujalak}) hold automatically.
\end{example}


\subsection{With inner migration}\label{withim}

Now imagine model (\ref{p:2}) with the inner migration of each patch, provided that the economy is closed. The spatial domain inside a country is one--dimensional again, as it was in case of the model (\ref{h:1rd}). The "length" of the $i$-th country is denoted by $L_i$, and $[0,L_i]\subset\mathbb{R}$, $L_i>0$, $i=1,2$. We assume that in each country, capital and labour are increased and decreased by the in- and outflowed quantity. Besides there is inner migration as well. (The situation is similar to the one which occurs often in hydrodynamics when a source and a "sink" are in a closed domain.) Moreover, we suppose that in case of migration an employee prefer such area in the other country which is situated the same distance from the border of the countries. More precisely, let us consider an employee living in place $x\in [0,L_1]$. If it will move to another country then it will choose exactly the place $y=x\frac{L_2}{L_1}\in [0,L_2]$ first, to find job.
Considering these ideas, the complete model can be written as:
\begin{equation}\label{pi:1hat}\left.\begin{array}{lll}
u'_t(1,t,x)=d_{11}(1) u''_{xx}(1,t,x)-d_{12}(1) v''_{xx}(1,t,x)+\tilde{f}_1(u(1,t,x),v(1,t,x),\hat{u}(2,t,x\frac{L_2}{L_1}),\hat{v}(2,t,x\frac{L_2}{L_1}))\\
v'_t(1,t,x)=-d_{21}(1) u''_{xx}(1,t,x)+d_{22}(1) v''_{xx}(1,t,x)+\tilde{g}_1(u(1,t,x),v(1,t,x),\hat{u}(2,t,x\frac{L_2}{L_1}),\hat{v}(2,t,\frac{L_2}{L_1}x))\\
\hat{u}'_t(2,t,y)=\hat{d}_{11}(2) \hat{u}''_{xx}(2,t,y)-\hat{d}_{12}(2) \hat{v}''_{xx}(2,t,y)+\tilde{f}_2(u(1,t,y\frac{L_1}{L_2}),v(1,t,y\frac{L_1}{L_2}),\hat{u}(2,t,y),\hat{v}(2,t,y))\\
\hat{v}'_t(2,t,y)=-\hat{d}_{21}(2) \hat{u}''_{xx}(2,t,y)+\hat{d}_{22}(2) \hat{v}''_{xx}(2,t,y)+\tilde{g}_2(u(1,t,y\frac{L_1}{L_2}),v(1,t,y\frac{L_1}{L_2}),\hat{u}(2,t,y),\hat{v}(2,t,y))
\end{array}\right\},\end{equation}
where
\begin{equation}\label{pi:1ahat}\begin{array}{l}
\tilde{f}_1(u(1,t,x),v(1,t,x),\hat{u}(2,t,x\frac{L_2}{L_1}),\hat{v}(2,t,x\frac{L_2}{L_1}))\\
=f(u(1,t,x),v(1,t,x))+\delta_1 (\rho_1(\hat{v}(2,t,x\frac{L_2}{L_1}))\hat{u}(2,t,x\frac{L_2}{L_1})-\rho_1(v(1,t,x))u(1,t,x))\\
\tilde{g}_1(u(1,t,x),v(1,t,x),\hat{u}(2,t,x\frac{L_2}{L_1}),\hat{v}(2,t,x\frac{L_2}{L_1}))\\
=g(u(1,t,x),v(1,t,x))+\delta_2 (\rho_2(\hat{u}(2,t,x\frac{L_2}{L_1}))\hat{v}(2,t,x\frac{L_2}{L_1})-\rho_2(u(1,t,x))v(1,t,x))\\
\tilde{f}_2(u(1,t,y\frac{L_1}{L_2}),v(1,t,y\frac{L_1}{L_2}),\hat{u}(2,t,y),\hat{v}(2,t,y))\\
=f(\hat{u}(2,t,y),\hat{v}(2,t,y))+\delta_1 (\rho_1(v(1,t,y\frac{L_1}{L_2}))u(1,t,y\frac{L_1}{L_2})-\rho_1(\hat{v}(2,t,y))\hat{u}(2,t,y))\\
\tilde{g}_2(u(1,t,y\frac{L_1}{L_2}),v(1,t,y\frac{L_1}{L_2}),\hat{u}(2,t,y),\hat{v}(2,t,y))\\
=g(\hat{u}(2,t,y),\hat{v}(2,t,y))+\delta_2 (\rho_2(u(1,t,y\frac{L_1}{L_2}))v(1,t,y\frac{L_1}{L_2})-\rho_2(\hat{u}(2,t,y))\hat{v}(2,t,y))
,\end{array}\end{equation}
where $f(u(1,t,x),v(1,t,x)))$ and $f(\hat{u}(2,t,y),\hat{v}(2,t,y)))$  are given by (\ref{p:3}) .
The  boundary conditions are similar to (\ref{h:1a}):
\begin{equation}\label{boundaryhat}\begin{array}{l}
u^{'}_x(1,t,0)=u^{'}_x(1,t,L_1)=0=v^{'}_x(1,t,0)=v^{'}_x(1,t,L_1),\\
\hat{u}^{'}_y(2,t,0)=\hat{u}^{'}_y(2,t,L_2)=0=\hat{v}^{'}_y(2,t,0)=\hat{v}^{'}_y(2,t,L_2).
\end{array}\end{equation}
 Let us change variable $y$ by $x\frac{L_2}{L_1}$ and introduce $\hat{u}(2,t,y)=\hat{u}(2,t,x\frac{L_2}{L_1})=u(2,t,x)$, $\hat{v}(2,t,y)=\hat{v}(2,t,x\frac{L_2}{L_1})=v(2,t,x)$, $d_{ij}(2)=\hat{d}_{ij}(2)(\frac{L_1}{L_2})^2$. Thus, we obtain:
 \begin{equation}\label{pi:1}\left.\begin{array}{lll}
u'_t(1,t,x)=d_{11}(1) u''_{xx}(1,t,x)-d_{12}(1) v''_{xx}(1,t,x)+\tilde{f}_1(u(1,t,x),v(1,t,x),u(2,t,x),v(2,t,x))\\
v'_t(1,t,x)=-d_{21}(1) u''_{xx}(1,t,x)+d_{22}(1) v''_{xx}(1,t,x)+\tilde{g}_1(u(1,t,x),v(1,t,x),u(2,t,x),v(2,t,x))\\
u'_t(2,t,x)=d_{11}(2) u''_{xx}(2,t,x)-d_{12}(2) v''_{xx}(2,t,x)+\tilde{f}_2(u(1,t,x),v(1,t,x),u(2,t,x),v(2,t,x))\\
v'_t(2,t,x)=-d_{21}(2) u''_{xx}(2,t,x)+d_{22}(2) v''_{xx}(2,t,x)+\tilde{g}_2(u(1,t,x),v(1,t,x),u(2,t,x),v(2,t,x))
\end{array}\right\},\end{equation}
where
\begin{equation}\label{pi:1a}\begin{array}{l}
\tilde{f}_1(u(1,t,x),v(1,t,x),u(2,t,x),v(2,t,x))\\
=f(u(1,t,x),v(1,t,x))+\delta_1 (\rho_1(v(2,t,x))u(2,t,x)-\rho_1(v(1,t,x))u(1,t,x))\\
\tilde{g}_1(u(1,t,x),v(1,t,x),u(2,t,x),v(2,t,x))\\
=g(u(1,t,x),v(1,t,x))+\delta_2 (\rho_2(u(2,t,x))v(2,t,x)-\rho_2(u(1,t,x))v(1,t,x))\\
\tilde{f}_2(u(1,t,x),v(1,t,x),u(2,t,x),v(2,t,x))\\
=f(u(2,t,x),v(2,t,x))+\delta_1 (\rho_1(v(1,t,x))u(1,t,x)-\rho_1(v(2,t,x))u(2,t,x))\\
\tilde{g}_2(u(1,t,x),v(1,t,x),u(2,t,x),v(2,t,x))\\
=g(u(2,t,x),v(2,t,x))+\delta_2 (\rho_2(u(1,t,x))v(1,t,x)-\rho_2(u(2,t,x))v(2,t,x))
,\end{array}\end{equation}
where $f(u(i,t,x),v(i,t,x)))$ are given by (\ref{p:3}) ($i=1,2$).
The  boundary conditions are similar to (\ref{h:1a}):
\begin{equation}\label{boundary}\begin{array}{l}
u^{'}_x(1,t,0)=u^{'}_x(1,t,L_1)=0=v^{'}_x(1,t,0)=v^{'}_x(1,t,L_1),\\
u^{'}_x(2,t,0)=u^{'}_x(2,t,L_1)=0=v^{'}_x(2,t,0)=v^{'}_x(2,t,L_1).
\end{array}\end{equation}
This latter is equivalant to
$\hat{u}^{'}_y(2,t,0)=\hat{u}^{'}_y(2,t,L_2)=0=\hat{v}^{'}_y(2,t,0)=\hat{v}^{'}_y(2,t,L_2)$.
The non--trivial equilibrium point $\bar{E_{p}}$ of system (\ref{p:2}) is a spatially constant solution of system (\ref{pi:1}) (and  of (\ref{pi:1hat}), as well) , because
\begin{equation*}
f(\bar{u},\bar{v})=0=g(\bar{u},\bar{v}),
\end{equation*}
(where $(\bar{u},\bar{v})$ is given by (\ref{eqrd})) and
substituting $\bar{E_p}=(\bar{u},\bar{v},\bar{u},\bar{v})$ into (\ref{pi:1a}):
\begin{equation*}\left.
\begin{array}{lll}
\tilde{f}_i(\bar{u},\bar{v},\bar{u},\bar{v})&=&\underbrace{f(\bar{u},\bar{v})}_{=0}+\delta_1 \underbrace{(\rho_1(\bar{v})\bar{u}-\rho_1(\bar{v})\bar{u})}_{=0}=0\\
\tilde{g}_i(\bar{u},\bar{v},\bar{u},\bar{v})&=&\underbrace{g(\bar{u},\bar{v})}_{=0}+\delta_2 \underbrace{(\rho_2(\bar{u})\bar{v}-\rho_2(\bar{u})\bar{v})}_{=0}=0
\end{array}\right\},\;\; i=1,2.\end{equation*}
\par Let we introduce the following notations: $U=\left[\begin{array}{cccc} u(1,t,x),&v(1,t,x),&u(2,t,x),&v(2,t,x)\end{array}\right]^{\mathbf{T}}$ and $F(U)=\left[\begin{array}{cccc} \tilde{f}_1,&\tilde{g}_1,&\tilde{f}_2,&\tilde{g}_2\end{array}\right]^{\mathbf{T}}$. Thus, we can transform model (\ref{pi:1}) to the following form:
\begin{equation}\label{pi:2}
U^{'}_t=D_4U^{''}_{xx}+F(U),
\end{equation}
where the diffusion matrix is: $\mathbf{D_4}=\left(\begin{array}{cccc}d_{11}(1)&-d_{12}(1)&0&0\\-d_{21}(1)&d_{22}(1)&0&0\\0&0&d_{11}(2)&-d_{12}(2)\\0&0&-d_{21}(2)&d_{22}(2)\end{array}\right).$
Linearizing system (\ref{pi:2}) at $\bar{E_p}$ and introducing the notation $V=U-\left[\begin{array}{cccc}\bar{u},&\bar{v},&\bar{u},&\bar{v}\end{array}\right]^{\mathbf{T}}$, we obtain:
\begin{equation}\label{pi:3}
V^{'}_t=\mathbf{D_4}V^{''}_{xx}+(\mathbf{A_p+\Gamma})V.
\end{equation}
We can solve equation (\ref{pi:3}) by Fourier--method again. It is easy to see that the stability of $\bar{E_p}$ depends on the stability properties of the matrix $(\mathbf{A_p}+\Gamma-\lambda_k\mathbf{D_4})$, where $\lambda_k=(\frac{k\pi}{L_1})^2$, $k=0,1,2,\dots$.
In order to get its characteristic polinomial in a much simplier form, we use the following row- and column operations:
\begin{equation*}\begin{split}
\det\left(\mathbf{A_p}+\Gamma-\lambda_k\mathbf{D}-\mu\mathbf{I}\right)=\\
\scriptsize\left|\begin{array}{cccc}
a_{r_{11}}-\delta_1\rho_1(\bar{v})-\lambda_kd_{11}(1)-\mu & a_{r_{12}}-\delta_1\rho'_1(\bar{v})\bar{u}+\lambda_kd_{12}(1) & \delta_1\rho_1(\bar{v}) & \delta_1\rho'_1(\bar{v})\bar{u}\\
a_{r_{21}}-\delta_2\rho'_2(\bar{u})\bar{v}+\lambda_kd_{21}(1) & a_{22}-\delta_2\rho_2(\bar{u})-\lambda_kd_{22}(1) -\mu & \delta_2\rho'_2(\bar{u})\bar{v} & \delta_2\rho_2(\bar{u})\\
\delta_1\rho_1(\bar{v}) & \delta_1\rho'_1(\bar{v})\bar{u} & a_{r_{11}}-\delta_1\rho_1(\bar{v})-\lambda_kd_{11}(2)-\mu & a_{r_{12}}-\delta_1\rho'_1(\bar{v})\bar{u}+\lambda_kd_{12}(2)\\
\delta_2\rho'_2(\bar{u})\bar{v} & \delta_2\rho_2(\bar{u}) & a_{r_{21}}-\delta_2\rho'_2(\bar{u})\bar{v}+\lambda_kd_{21}(2) & a_{r_{22}}-\delta_2\rho_2(\bar{u})-\lambda_kd_{22}(2)-\mu\\
\end{array}\right|\\
\scriptsize\underbrace{=}_{\begin{array}{c}c1+c3;\\c2+c4\\r3-r1;\\r4-r2\end{array}}\scriptsize\left|\begin{array}{cccc}
a_{r_{11}}-\lambda_kd_{11}(1)-\mu & a_{r_{12}}+\lambda_kd_{12}(1) & \delta_1\rho_1(\bar{v}) & \delta_1\rho'_1(\bar{v})\bar{u}\\
a_{r_{21}}+\lambda_kd_{21}(1) & a_{r_{22}}-\lambda_kd_{22}(1)-\mu & \delta_2\rho'_2(\bar{u})\bar{v} & \delta_2\rho_2(\bar{u})\\
-\lambda_k(d_{11}(2)-d_{11}(1)) & -\lambda_k(d_{12}(2)-d_{12}(1)) & a_{r_{11}}-2\delta_1\rho_1(\bar{v})-\lambda_kd_{11}(2)-\mu & a_{12}-2\delta_1\rho'_1(\bar{v})\bar{u}+\lambda_kd_{12}(2)\\
-\lambda_k(d_{21}(2)-d_{21}(1))& -\lambda_k(d_{22}(2)-d_{22}(1)) & a_{r_{21}}-2\delta_2\rho'_2(\bar{u})\bar{v}+\lambda_kd_{21}(2) & a_{r_{22}}-2\delta_2\rho_2(\bar{u})-\lambda_kd_{22}(2)-\mu\\
\end{array}\right|.
\end{split}\end{equation*}
Let us consider the special case when there is no difference between the velocities of the diffusion--coeffitients inside the two countries (patches). This is valid for the self-- and the cross--diffusion. Thus,
\begin{equation}\label{d-k}
d_{ij}(2)=d_{ij}(1):=d_{ij}, \;\;i,j=1,2.
\end{equation}
In this case the latter determinant is equal to:
\begin{equation}
\scriptsize
\det\left[\begin{array}{cc}
a_{r_{11}}-\lambda_kd_{11}-\mu & a_{r_{12}}+\lambda_kd_{12}\\
a_{r_{21}}+\lambda_kd_{21} & a_{r_{22}}-\lambda_kd_{22}-\mu
\end{array}\right]
\det\left[\begin{array}{cc}
a_{r_{11}}-2\delta_1\rho_1(\bar{v})-\lambda_kd_{11}-\mu & a_{r_{12}}-2\delta_1\rho'_1(\bar{v})\bar{u}+\lambda_kd_{12}\\
{a_{r_{21}}-2\delta_2\rho'_2(\bar{u})\bar{v}+\lambda_kd_{21}} & a_{r_{22}}-2\delta_2\rho_2(\bar{u})-\lambda_kd_{22}-\mu\end{array}\right],
\end{equation}
Taking into account the signs of the coefficients: $d_{ij}>0$, $\rho_i>0$, $\delta_i>0$, $\rho'_1(\bar{v})<0$, $\rho'_2(\bar{u})<0$, and assuming (\ref{plusmas}) we obtain:
\begin{equation*}
\scriptsize
\det\left[\begin{array}{cc}
\underbrace{a_{r_{11}}-\lambda_kd_{11}}_{<0}-\mu & a_{r_{12}}+\lambda_kd_{12}\\
\underbrace{a_{r_{21}}+\lambda_kd_{21}}_{>0} & \underbrace{a_{r_{22}}-\lambda_kd_{22}}_{<0}-\mu
\end{array}\right]
\det\left[\begin{array}{cc}
\underbrace{a_{r_{11}}-2\delta_1\rho_1(\bar{v})-\lambda_kd_{11}}_{<0}-\mu & a_{12}-2\delta_1\rho'_1(\bar{v})\bar{u}+\lambda_kd_{12}\\
\underbrace{a_{r_{21}}-2\delta_2\rho'_2(\bar{u})\bar{v}+\lambda_kd_{21}}_{>0} & \underbrace{a_{r_{22}}-2\delta_2\rho_2(\bar{u})-\lambda_kd_{22}}_{<0}-\mu\end{array}\right].
\end{equation*}
In order to get a stable polinomial, these two determinants have to be separately stable polinomials depending on $\mu$. Since the entries in the main diagonal are negative, the two determinants (without subtracting $\mu$ in the main diagonals) have to be positive. Namely the following two conditions have to be satisfied:
\begin{equation}\label{det1}
\begin{split}
(a_{r_{11}}-\lambda_kd_{11})(a_{r_{22}}-\lambda_kd_{22})-(a_{r_{12}}+\lambda_kd_{12})(a_{r_{21}}+\lambda_kd_{21})\\
 =\lambda_k^2\det D-\lambda_k(a_{r_{11}} d_{22}+a_{r_{22}} d_{11}+a_{r_{12}} d_{21}+a_{r_{21}} d_{12})+\det A_r>0.
\end{split}
\end{equation}

\begin{equation}\label{det2}\begin{split}
(a_{r_{11}}-2\delta_1\rho_1(\bar{v})-\lambda_kd_{11})(a_{r_{22}}-2\delta_2\rho_2(\bar{u})-\lambda_kd_{22})\\
-(a_{12}-2\delta_1\rho'_1(\bar{v})\bar{u}+\lambda_kd_{12})(a_{r_{21}}-2\delta_2\rho'_2(\bar{u})\bar{v}+\lambda_kd_{21})\\
=\lambda_k^2 \det D\\
-\lambda_k [(a_{r_{11}} d_{22}+a_{r_{22}} d_{11}+a_{r_{12}} d_{21}+a_{r_{21}}d_{12})]\\
+2[\delta_2\rho_2(\bar{u})d_{11}+\delta_1\rho_1(\bar{v})d_{22}+\delta_2\rho'_2(\bar{u})\bar{v}d_{12}+\delta_1\rho'_1(\bar{v})\bar{u}d_{21}])\\
+(\det A_r-2(\delta_1\rho_1(\bar{v})a_{r_{22}}+\delta_2\rho_2(\bar{u})a_{r_{11}}-\delta_1\rho'_1(\bar{v})\bar{u}a_{r_{21}}-\delta_2\rho'_2(\bar{u})\bar{v}a_{r_{12}})\\
+4\delta_1\delta_2(\rho_1(\bar{v})\rho_2(\bar{u})-\rho'_1(\bar{v})\bar{u}\rho'_2(\bar{u})\bar{v}))>0.
\end{split}
\end{equation}
In connection with (\ref{det1}) it can be seen that if (\ref{h:7rd}) is satisfied then it is a stable polynomial, thus it is positive for all $\lambda_k\geq0$. \\
In connection with (\ref{det2}) it can be seen that if (\ref{h:7rd}), (\ref{1.feltetel}), (\ref{2.feltetelujalak}) and
\begin{eqnarray}\label{pathdkeplet}
\delta_2\rho_2(\bar{u})d_{11}+\delta_1\rho_1(\bar{v})d_{22}+\delta_2\rho'_2(\bar{u})\bar{v}d_{12}+\delta_1\rho'_1(\bar{v})\bar{u}d_{21}>0
\end{eqnarray}
are satisfied then it is also a stable polynomial, thus it is positive for all $\lambda_k\geq0$ again. \\
Condition (\ref{pathdkeplet}) can be rewritten in the following form:
\begin{eqnarray}\label{pathdkepletmas}
d_{12}<\frac{\delta_2\rho_2(\bar{u})d_{11}+\delta_1\rho_1(\bar{v})d_{22}+\delta_1\rho'_1(\bar{v})\bar{u}d_{21}}{-\delta_2\rho'_2(\bar{u})\bar{v}}.
\end{eqnarray}

We have arrived at:
\begin{thm}\label{TH4.4}
Suppose that conditions (\ref{h:3rd}), (\ref{h:2rd}), (\ref{h:4rd}), (\ref{h:7rd}), (\ref{plusmas}), (\ref{1.feltetel}), (\ref{2.feltetelujalak}), (\ref{d-k}), (\ref{pathdkepletmas}), $\rho'_{i}<0, \;\;\rho_{i}>0, \;\; i=1,2$ hold for each patch then
equilibrium $\bar{E_p}=(\bar{u},\bar{v},\bar{u},\bar{v})$  of system (\ref{pi:1hat}) (given by (\ref{eqrd})) is asymptotically stable.
\end{thm}
This theorem means that if we consider a country as a patch it is required not to have a too low number of free jobs in a country (according to (\ref{plusmas})); the hectic flow of the capital towards places where labour force is dense has to be avoided inside the country (according to (\ref{h:7rd}) and (\ref{pathdkepletmas})); the migration functions $\rho_i$, $i=1,2$ has to be characterized by ((\ref{1.feltetel}), (\ref{2.feltetelujalak})) that describe the inflow of capital and labour into the two countries (patches), and those give a limit for inflow of capital, the two countries have similar market practices (according to (\ref{d-k})) then the common market is stable.

Comparing the conditions of Theorem {\ref{TH4.3}} and Theorem {\ref{TH4.4}} it can be seen that it is not enough to guarantee the stability of the markets between the two countries. If we consider the spatial movements inside the countries then stability conditions are required more stronger. In other words, not properly designed or supervised trade mechanism between the respective countries may as well cause the destabilization of each market, irrespectively of their own stability. Now the reversed case might also happen, namely due to a change in any (or both) market designs, the originally balanced (stable) joint market can lose its stability. All this means, that internal (external) market changes cannot occur without the distortion of the external (internal) market.

\begin{rem}
By \emph{Casten--Holland (1977)} \cite{Casten Holland}, the asymptotic stability of the constant solution $(\bar{u},\bar{v},\bar{u},\bar{v})$ of system (\ref{pi:1hat}) is implied.
\end{rem}
Increasing diffusion coefficient $d_{12}$  system (\ref{pi:1hat}) may lose its stability.

\section{Conclusion}\label{sum}

 In section \ref{ratio} the stability of the market in a country was considered in the case of diffusion of capital and labour. It was established that the hectic movement of the capital stock toward places where labour force is dense acts negatively and able to destabilize the system. We determined a critical value for the cross diffusion coefficient ($d_{12}$) under which the market is stable in the country. Then  in section \ref{patch} two similar country was considered with their own stable markets. We studied the effect of the movement of capital and labour between the countries to the stability behaviour. In this case the hectic movement of the capital inside a country has to avoid even more. It may be necessary to reduce the inner flow of capital toward labour more significantly. This result is not surprising, since capital flow, namely a new investment, in order to obtain new jobs has to be slower than the movement of the labour.
  As a conclusion, the common market can be unstable if the inner movement of capital is sufficiently fast. This movement has the same effect on the individual market of each countrie, hence the destabilization of each countrie may affect the destabilization of the common market, and vice versa: instability of the common market makes the individual markets to be unstable.
The wise moderation of capital is useful in all cases.
Our results are applicable in case of two similar countries, in those where the velocity of capital flow and labour movement is approximately equal in the two countries, namely where there is similar market behaviour. Thus, our investigation can be applied for EU members, but not for an EU member and China common market. In this latter case we have to reject condition (\ref{d-k}) because the market practices are different and the stability investigation is become more complicated.

\end{document}